\documentclass{gtmon_a}
\pdfoutput=1

\proceedingstitle{Proceedings of the School and Conference in Algebraic
Topology (The Vietnam National University, Hanoi, 9-20 August 2004)}
\conferencestart{9 August 2004}
\conferenceend{20 August 2004}
\conferencename{School and Conference in Algebraic Topology}
\conferencelocation{Vietnam National University, Hanoi, Vietnam}

\editor{John Hubbuck}
\givenname{John}
\surname{Hubbuck}

\editor{Nguy\~\ecircumflex{}n H V H\uhorn{}ng}
\givenname{H\uhorn{}ng}
\surname{Nguy\~\ecircumflex{}n}

\editor{Lionel Schwartz}
\givenname{Lionel}
\surname{Schwartz}

\title{Young tableaux and the Steenrod algebra}

\author{G\,Walker}
\givenname{Grant}
\surname{Walker}
\address{School of Mathematics\\
University of Manchester\\\newline
Manchester\\
M13 9PL\\
United Kingdom}
\email{Grant.Walker@manchester.ac.uk}
\urladdr{}

\author{R\,M\,W Wood}
\givenname{R M W}
\surname{Wood}
\email{Reg.Wood@manchester.ac.uk}
\urladdr{}

\volumenumber{11}
\issuenumber{}
\publicationyear{2007}
\papernumber{18}
\startpage{379}
\endpage{397}

\doi{}
\MR{}
\Zbl{}

\arxivreference{}



\received{28 December 2004}
\revised{17 July 2005}
\accepted{18 December 2005}
\published{14 November 2007}
\publishedonline{14 November 2007}
\proposed{}
\seconded{}
\corresponding{}
\editor{}
\version{}

\newcommand{\Sq}{\mathrm{Sq}}
\newcommand{\cp}{\mathrm{cp}}
\renewcommand{\theenumi}{\roman{enumi}}

\makeatletter
\def\cnewtheorem#1[#2]#3{\newtheorem{#1}{#3}[section]
\expandafter\let\csname c@#1\endcsname\c@theorem}

\newtheorem{theorem}{Theorem}[section]
\cnewtheorem{corollary}[theorem]{Corollary}
\cnewtheorem{proposition}[theorem]{Proposition}
\cnewtheorem{lemma}[theorem]{Lemma}
\cnewtheorem{conjecture}[theorem]{Conjecture}
\theoremstyle{remark}
\cnewtheorem{problem}[theorem]{Problem}
\cnewtheorem{problems}[theorem]{Problems}
\cnewtheorem{example}[theorem]{Example}
\cnewtheorem{remark}[theorem]{Remark}
\cnewtheorem{definition}[theorem]{Definition}

\makeatother  

\begin{document}

\begin{asciiabstract} 
The purpose of this paper is to forge a direct link between
the hit problem for the action of the Steenrod algebra A  on the
polynomial algebra P(n)=F_2[x_1,...,x_n],  over the field F_2 of two
elements,  and semistandard Young tableaux as they apply to the modular
representation theory of the  general linear group GL(n,F_2). The cohits
Q^d(n)=P^d(n)/P^d(n)\cap A^+(P(n)) form a modular representation of
GL(n,F_2) and the hit problem is to analyze this module.   In certain
generic degrees d we show how the semistandard Young tableaux can be
used to index a set of monomials which span Q^d(n).  The hook formula,
which calculates the number of semistandard  Young tableaux, then gives an
upper bound for the dimension of Q^d(n).  In the particular degree d where
the Steinberg module appears for the first time in P(n) the upper bound
is exact and Q^d(n) can then be identified with the Steinberg module.
\end{asciiabstract}

\begin{htmlabstract}
The purpose of this paper is to forge a direct link between the hit
problem for the action of the Steenrod algebra A on the
polynomial algebra <b>P</b>(n)=<b>F</b><sub>2</sub>[x<sub>1</sub>,&hellip;,x<sub>n</sub>],  over the
field  <b>F</b><sub>2</sub> of two elements,  and semistandard Young tableaux
as they apply to the modular representation theory of the  general linear
group GL(n,<b>F</b><sub>2</sub>). The cohits
<b>Q</b><sup>d</sup>(n)=<b>P</b><sup>d</sup>(n)/<b>P</b><sup>d</sup>(n)&cap;A<sup>+</sup>(<b>P</b>(n)) form a modular representation of
GL(n,<b>F</b><sub>2</sub>) and the hit problem is to analyze this module.
In certain generic degrees d we show how the semistandard Young tableaux
can be used to index a set of monomials which span  <b>Q</b><sup>d</sup>(n).
The hook formula, which calculates the number of semistandard  Young
tableaux, then gives an upper bound for the dimension of  <b>Q</b><sup>d</sup>(n).
In the particular degree d where the Steinberg module appears for the
first time in <b>P</b>(n) the upper bound is exact and  <b>Q</b><sup>d</sup>(n)
can then be identified with the Steinberg module.
\end{htmlabstract}

\begin{abstract} 
The purpose of this paper is to forge a direct link between the hit
problem for the action of the Steenrod algebra ${\cal A}$  on the
polynomial algebra ${\bf P}(n)={\mathbb F}_2[x_1,\ldots,x_n]$,  over the
field $ {\mathbb F}_2$ of two elements,  and semistandard Young tableaux
as they apply to the modular representation theory of the  general linear
group $GL(n,{\mathbb F}_2)$. The cohits ${\bf Q}^d(n)={\bf P}^d(n)/{\bf
P}^d(n)\cap{\cal A}^+({\bf P}(n))$ form a modular representation of
$GL(n,{\mathbb F}_2)$ and the hit problem is to analyze this module.
In certain generic degrees $d$ we show how the semistandard Young tableaux
can be used to index a set of monomials which span  ${\bf Q}^d(n)$.
The hook formula, which calculates the number of semistandard  Young
tableaux, then gives an upper bound for the dimension of  ${\bf Q}^d(n)$.
In the particular degree $d$ where the Steinberg module appears for the
first time in ${\bf P}(n)$ the upper bound is exact and  ${\bf Q}^d(n)$
can then be identified with the Steinberg module.
\end{abstract}

\maketitle


\section{Introduction}
\label{sec:sec1}

Young tableaux form  a combinatorial device for constructing 
representations   of the general linear group  $GL(n)$ of $n\times n$ 
non-singular matrices and its subgroup  $\Sigma_n$ of permutation 
matrices,  both in the classical  case,  over the field of complex 
numbers,  and in  the modular case,  where the characteristic of the 
field divides the order of the group (see Fulton \cite{f},
James--Kerber \cite{j}, MacDonald \cite{m}, Sagan \cite{s1}
and Stanley \cite{s2}). The group  
$GL(n,{\mathbb F}_2)$, over the field $ {\mathbb F}_2$ of two 
elements, acts naturally on  the polynomial  algebra 
${\bf P}(n)={\mathbb F}_2[x_1,\ldots,x_n]$ 
by matrix substitution and the 
homogeneous polynomials   ${\bf P}^d(n)$ of degree $d$ form  a 
representation space.  The modular
representation theory  of subgroups of  $GL(n,{\mathbb F}_2)$, acting 
in this way on  ${\bf P}(n)$,  is 
important in understanding  the nature of the \emph{hit problem\/} for 
the  action of the mod $2$  Steenrod algebra  ${\cal A}$ on 
${\bf P}(n)$.  The problem  is to find a minimal 
generating set for ${\bf P}(n)$ as an  ${\cal A}$--module
(see Boardman \cite{b},
Janfada--Wood \cite{jw1,jw2}, Kameko \cite{k1,k2}, Peterson \cite{p} and
Wood \cite{w1,w2,w3,w4,w5}). The Steenrod squaring operators
$\Sq^k$ generate ${\cal A}$ as an algebra and act as  
$GL(n,{\mathbb F}_2)$--module maps from   ${\bf P}^d(n)$ to  
${\bf P}^{d+k}(n)$. A polynomial $h$ is \emph{hit\/} if it can be 
written as a finite sum  $h=\sum_{k>0}\Theta_k(f_k)$ for elements 
$\Theta_k$   of positive grading in ${\cal A}$ and suitable 
polynomials $f_k$. Equivalently,  $h=\sum_{k>0}\Sq^k(g_k)$ for 
suitable polynomials $g_k$.

The hit problem can be viewed in terms of finding a vector space 
basis for the quotient ${\bf Q}(n)$  of ${\bf P}(n)$ by the hit 
elements. This quotient is a $GL(n,{\mathbb F}_2)$--module.  It is 
clear that polynomials in  ${\bf P}^d(n)$ which represent 
non-trivial elements in an irreducible  composition factor of   
${\bf P}^d(n)$, occurring for the first time in degree $d$,  cannot be hit, 
otherwise there would be a Steenrod operation linking the composition 
factor with an earlier occurrence. This goes some way to explain the 
interrelationship between
modular representation theory and the Steenrod algebra. From the 
point of view 
of representation theory, the quotient ${\bf Q}(n)$  is a 
repository for the irreducible modular representations  of   
$GL(n,{\mathbb F}_2)$ and from the point of view of the Steenrod 
algebra, first occurrences of irreducible representations contribute 
to a generating set of the ${\cal A}$--module ${\bf P}(n)$.

We shall explain  a direct connection between 
Young tableaux and generators for the  ${\cal A}$--module  
${\bf P}(n)$. In general, there are  too many Young tableaux  to solve the 
hit problem
precisely  but in certain degrees the number of semistandard Young 
tableaux  does  give the correct minimal  number of generators. The 
following is a sample result from 
the  more  general \fullref{main}. 

\begin{theorem}\label{th01}
In any minimal generating set for the  
${\cal A}$--module  ${\bf P}(n)$, there are 
$2^{\binom{n}{2}}$ elements in degree $d=2^n -n -1$. In this degree
monomial  generators  in  ${\bf P}(n)$ may be chosen in bijective 
correspondence with the semistandard Young tableaux 
associated with the partition $(n-1,n-2,\ldots,1)$ of the number
$\binom{n}{2}$. Furthermore,  these generators provide 
representatives for an additive basis for the first occurrence  in 
degree $d$   of  the Steinberg representation of  
$GL(n,{\mathbb F}_2)$, viewed  as the  quotient of  ${\bf P}^d(n)$ by the hit 
elements. 
\end{theorem}
For instance, taking the case $n=3$, there are eight 
 semistandard Young tableaux as exhibited below.
\begin{example}\label{ex02}
$$
\begin{smallmatrix}
1&1\\
2
\end{smallmatrix}\qua 
\begin{smallmatrix}
1&1\\
3
\end{smallmatrix}\qua
\begin{smallmatrix}
1&2\\
2
\end{smallmatrix}\qua 
\begin{smallmatrix}
1&2\\
3
\end{smallmatrix}\qua 
\begin{smallmatrix}
2&2\\
3
\end{smallmatrix}\qua 
\begin{smallmatrix}
2&3\\
3
\end{smallmatrix}\qua 
\begin{smallmatrix}
1&3\\
2
\end{smallmatrix}\qua 
\begin{smallmatrix}
1&3\\
3
\end{smallmatrix}
$$
The  corresponding  monomial  generators  in   ${\bf P}^4(3)$, 
equivalently representative monomials of a vector space basis for  
${\bf Q}^4(3)$,  are 
$$x_1^3x_2,\qua x_1^3x_3,\qua x_1x_2^3,\qua
x_1x_2^2x_3,\qua x_2^3x_3,\qua x_2x_3^3,\qua 
x_1x_2x_3^2,\qua x_1x_3^3.$$
\end{example}
The fact that the first occurrence of the Steinberg representation is 
in degree $2^n-n-1$ is well known (Mitchell--Priddy \cite{mp},
Minh--Tri \cite{mt} and Walker--Wood \cite{ww}). The result of 
\fullref{th01} may be paraphrased by saying that the Steinberg 
representation is the only irreducible representation of  
$GL(n,{\mathbb F}_2)$  to contribute to a minimal generating set for 
the  ${\cal A}$--module  ${\bf P}(n)$ in this degree.

In the next section we explain  how  to associate monomials with 
tableaux and, more generally,  we translate some of the traditional  
language used in the combinatorial theory of  tableaux into the 
language of \emph{block technology\/}, which is appropriate for handling 
the action of the Steenrod algebra. In particular, we  introduce a 
combinatorial procedure, called \emph{splicing\/},  which is used to 
replace a block by a formal sum of \emph{semistandard blocks\/}. This  
is analogous to the \emph{straightening\/} process  for bringing Young 
tableaux into standard form in the context of group rings, see Fulton
\cite{f}. 
In \fullref{sec:sec3} it is shown how splicing can be realized by the action 
of the Steenrod algebra and \fullref{main} is proved. 
In \fullref{sec:sec4} we show how \fullref{th01} follows by considering 
the special case of the Steinberg representation, using the hook 
formula to count the number of  semistandard Young tableaux. 

In general, the hook formula shows that for a fixed $n\ge 2$ and 
increasing $d$ the number of semistandard Young tableaux increases, 
whereas the dimension of  ${\bf Q}^d(n)$ is known to be bounded in 
$d$ for a given $n$, see Carlisle--Wood
\cite{cw}. It would be interesting to find a more 
restrictive condition on semistandard Young tableaux
which cuts  down the number of generators of  ${\bf Q}^d(n)$, 
at least in the
row-regular case,  to a number bounded in $d$ which  estimates  more 
closely the dimension of the cohits. It would also be interesting to 
investigate the dual hit problem and identify  a basis for the kernel 
of the down  Steenrod action in terms of the combinatorics of Young 
tableaux and the relationship with the ring of lines as  described in 
Alghamdi--Crabb--Hubbuck \cite{ach} and Crabb--Hubbuck
\cite{ch}. At the end of \fullref{sec:sec3} we give an example to show 
the limitations of the main theorem. In the last section we explain 
briefly 
how \fullref{th01} can be extended  to  other 
irreducible representations of  $GL(n,{\mathbb F}_2)$ having a 
certain 
affinity to the Steinberg representation. 

\section{Binary blocks and Young tableaux} 
\label{sec:sec2}
There are two frequently used numerical functions in the context of 
the hit problem. One is  the 
$\alpha$--function  $\alpha(d)$ of a positive  integer $d$,
which counts the  number of digits $1$ in the binary expansion of 
$d$,  and the other is the  
$\mu$--function  $\mu(d)$, which is the smallest number $k$ for which 
$d$ can be partitioned in exponential form  
$d=\sum_{i=1}^{k}(2^{\lambda_i}-1)$. We extend the definitions to 
cover $\alpha(0)=\mu(0)=0$.  In general, the 
exponential partition of a number $d$, with a given value of   
$\mu(d)$,  is not unique. For example $\mu(17)=3$ and 
$17=15+1+1=7+7+3$.

We are concerned with two types of partitions of  numbers: the 
exponential partition of $d$ as in the definition of the  
$\mu$--function and the ordinary partition  
$\lambda=(\lambda_1,\ldots, \lambda_n )$ of the number
$|\lambda|=\lambda_1+ \cdots +\lambda_n $, where 
$\lambda_1\ge\lambda_2 \ge \ldots \ge \lambda_n\ge 0$. The 
\emph{length\/} of the partition is the number of non-zero parts $\lambda_i$. 
In this article we 
reserve $n$ for the number of variables in  ${\bf P}^d(n)$ and  
restrict attention to  partitions of  length not greater than $n$.  
In combinatorics it is customary to illustrate the  partition  
$\lambda$  by a \emph{Ferrers diagram\/}, which is an array of boxes  
in  echelon shape with $\lambda_i$ boxes  in the $i$th row. The 
positions of the boxes are the \emph{nodes\/} of the Ferrers diagram.  A 
\emph{Young tableau\/} is a   Ferrers diagram in which each box is 
filled with a positive integer.  In particular, filling each box with 
the digit $1$ produces an array of the form  
$$F=\begin{smallmatrix}
1&\ldots&1&\ldots&1&\ldots&1\\
1&\ldots&1&\ldots&1&&\\
&\ldots&&&&&\\ 
1&\ldots&1&&&&
\end{smallmatrix}$$
with $\lambda_i$ contiguous digits $1$ in the $i$th row. We shall 
call this a \emph{Ferrers block\/} and  interpret it  in terms of  the 
exponential partition  
$d=\sum_{i=1}^{k}(2^{\lambda_i}-1)$, where now the rows of $F$ are 
the reverse binary expansions of the numbers $2^{\lambda_i}-1$ as 
read from left to right.   
More generally, a  \emph{binary block\/} is  a $(0,1)$--array  
associated with 
a monomial $f=x_1^{d_1}\ldots x_n^{d_n}$, whose entries   are  the 
digits, in  reversed binary expansion,  of the exponents $d_i$. 
Blocks were introduced 
in Carlisle--Wood
\cite{cw} as a graphical  device for keeping track of the action 
of Steenrod squares on  monomials and have been used in several 
places to exhibit minimal sets of  monomial generators 
(see Janfada--Wood \cite{jw2}).   
A formal sum of blocks corresponds to a polynomial over 
${\mathbb F}_2$ 
(ordinary addition of matrices is not used in this article). If 
we are working in ${\bf P}(n)$ then the number of rows in a block is 
$n$. In particular a Ferrers block may have zero rows at the bottom.  On the other hand the number $c$ of columns in  a block is 
not determined by the corresponding monomial.  We  adopt the 
convention of regarding two row-vectors of nonnegative integers as 
equivalent if they differ by trailing zeros and we  omit trailing 
zeros when convenient.  In particular the empty vector is identified 
with a vector of zero entries. The convention is extended  to arrays 
except that an empty row in indicated by a leading $0$. This is 
necessary to keep track of the positions of  missing variables in a 
monomial and maintain the number of rows  at $n$.  Under these 
conventions we have a bijective correspondence between blocks with 
$n$ rows and monomials in  ${\bf P}(n)$. 
The blocks associated with the monomials  in \fullref{ex02} are 
given in the following list.
\begin{example}\label{ex11}
$$
\begin{smallmatrix}
1&1\\
1\\
0
\end{smallmatrix}\qua
\begin{smallmatrix}
1&1\\
0\\
1
\end{smallmatrix}\qua
\begin{smallmatrix}
1\\
1&1\\
0
\end{smallmatrix}\qua 
\begin{smallmatrix}
1\\
0&1\\
1
\end{smallmatrix}\qua 
\begin{smallmatrix}
0\\
1&1\\
1
\end{smallmatrix}\qua 
\begin{smallmatrix}
0\\
1\\
1&1
\end{smallmatrix}\qua 
\begin{smallmatrix}
1\\
1\\
0&1
\end{smallmatrix}\qua 
\begin{smallmatrix}
1\\
0\\
1&1
\end{smallmatrix}
$$
\end{example}

In the context of the hit problem,  a  \emph{spike\/} in ${\bf P}(n)$ is 
a monomial of the form \linebreak 
$x_1^{d_1}x_2^{d_2}\cdots x_n^{d_n}$, where each $d_i$ has the form $
2^{\lambda_i} -1$. The corresponding block is a row permutation of 
the Ferrers block, with appropriate $0$--rows inserted.

A more compact way of designating a $(0,1)$--array  $F$  is to  form 
the  
corresponding    array $Y$ of nonnegative   integers by  the 
following rule. 
For each number  $j$ let $L$ denote the list  of row positions 
occupied by a digit $1$ in the $j$th column of   $F$, counting from 
the top row down. Then the $i$th  element of $L$ occupies position 
$(i,j)$ in $Y$. If a column of $F$ has no digits $1$ then the 
corresponding column of $Y$ has  zero  entries, keeping in mind the 
trailing zero convention for rows of the array.  By construction, the 
non-zero entries in a  column  of $Y$ are strictly increasing. 
\begin{definition}\label{cost}
An array is \emph{column strict\/} if the non-zero  entries of any 
column are strictly increasing from top down. The process of 
assigning the column-strict array $Y$  to the block $B$ is called the 
\emph{column-position correspondence\/} and is denoted by $Y=\cp(F)$.  
\end{definition}
Given a column-strict array $Y$ with no entry larger than $n$, then 
it is clear how to  constitute the block $F$ with $n$ rows  so that 
$\cp(F)=Y$. The column-position correspondence is therefore bijective 
between blocks and 
column-strict  arrays.  It is easy to see  that the blocks of  
\fullref{ex11} and the arrays of \fullref{ex02} are related by the 
$\cp$ correspondence  and this in turn  establishes the correspondence 
with  the  monomials in \fullref{ex02}.

We shall now translate 
some of the traditional language  of Young tableaux (see Fulton \cite{f},
Macdonald \cite{m}, Sagan \cite{s1} and Stanley \cite{s2})
into the language of block technology (see Janfada--Wood \cite{jw2}). 
The 
\emph{$\omega$--vector\/} of a block $F$ is the vector
$\omega(F)=(\omega_1,\ldots, \omega_c)$  of column sums of $F$. The 
\emph{$\alpha$--vector\/} of $F$ is the vector  
$\alpha(F)=(\alpha_1,\ldots, \alpha_n)$  of row sums. In the case of 
the Ferrers block associated with the partition 
$\lambda$ we have $\alpha(F)=\lambda$ and $\omega(F)=\lambda'$, the 
conjugate  of $\lambda$.  The \emph{degree\/} $d$ of a block $F$,
or associated  array $\cp(F)$,   means the degree of the corresponding
 monomial,  and this is a function of the $\omega$--vector given by 
$d=\sum_{j>0} \omega_j2^{j-1}$. In terms of  $\cp(F)$, the  $j$th 
entry $\omega_j$ of  the  $\omega$--vector is  the
number of  non-zero entries in the  $j$th column of  $\cp(F)$.
 
Of particular interest in this article are the monomials with 
\emph{descending\/}
$\omega$--vectors, meaning that  
$\omega_j\ge \omega_{j+1}$ for $j\ge 1$, keeping in mind the trailing 
zeros convention for vectors.   
All the blocks in \fullref{ex11} are of this type with 
$\omega$--vector $(2,1)$.   If a block $F$ has a descending 
$\omega$--vector then it is easy to see that the corresponding 
column-position array   $Y=\cp(F)$ is a column-strict Young tableau. 
In some parts of the literature 
column-strict is included in the definition of a Young tableau.    
One can easily check that the $i$th entry  of  $\alpha(F)$ is the 
number of repetitions of $i$ in $Y$. 

It follows that a monomial with descending  $\omega$--vector has a 
uniquely 
associated col\-umn-strict Young tableau via the column-position  
correspondence.  In combinatorics a  column-strict  Young tableau is 
called \emph{semistandard\/} if the rows are weakly  increasing. The 
Young tableaux in \fullref{ex02} are semistandard. 

The  following lemma, which is straightforward  to prove,  summarizes 
the situation so far. 
\begin{lemma}\label{le13}
Working in ${\bf P}(n)$, the column-position correspondence sets up a 
bijection between  monomials  with descending  
$\omega$--vectors     and column-strict  Young tableaux, with entries 
taken from the set $\{1,\ldots,n\}$, based on Ferrers blocks  with  
$n$ rows.  A  semistandard tableau $\cp(F)$ corresponds   to a block 
$F$  with the property that
$\omega(F[i])$ is descending for each $i$ in the range 
$1\le i \le n$, where  $F[i]$  denotes   the block  formed by taking the first 
$i$ rows of $F$. 
\end{lemma}
In the light of this lemma it is appropriate to make the following 
definition. 

\begin{definition}
A block $F$ with $n$  rows is \emph{semistandard\/} if  $\omega(F[i])$ 
is descending for each  subblock of  $F[i]$ for $1\le i \le n$. 
\end{definition}

The ultimate aim of this article  is to find a generating set for the 
$\cal A$--module  ${\bf P}(n)$ among semistandard blocks at least in 
certain degrees. We shall call a degree $d$ \emph{row-regular for $n$\/} 
if it has  an exponential partition  
$d=\sum_{i=1}^{n}(2^{\lambda_i}-1)$, where the partition 
$\lambda=(\lambda_1,\ldots,  \lambda_{n})$  satisfies the condition   
$\lambda_1 > \cdots > \lambda_{n}\ge 0$. In this case $d$ has a 
unique  exponential partition of length $n-1$ or $n$. For such a 
degree we have $\mu(d)=n$ or 
$\mu(d)=n-1$,  but not all degrees with these $\mu$--values are 
row-regular for $n$. Up to permutation of rows there is only one 
spike  in  ${\bf P}^d(n)$ when $d$ is row-regular,   and therefore 
just one  associated Ferrers block  and just  one descending 
$\omega$--vector.   All monomials in  ${\bf P}^d(n)$ with this 
$\omega$--vector have corresponding column-strict  Young tableaux with 
the same   underlying Ferrers diagram. 

Later proofs will require induction on certain partial order 
relations  on monomials and corresponding blocks. These are  
constructed from total order relations on  $\omega$--vectors. We shall 
highlight two of these.
\begin{definition} 
Let $L=(a_1,a_2,\ldots,a_s)$ and $M=(b_1,b_2,\ldots,b_s)$ be two 
vectors
of non-negative integers.  We write $L>_lM$ and read  `greater than 
in left order'  if $a_1>b_1$ or $a_i=b_i$ for $1\le i < t\le s$ and 
$a_t>b_t$.  We also write  $L>_r M$ and read `greater than  in  right 
order'  if $a_s<b_s$ or $a_i=b_i$ for 
$1\le t<i\le s$ and $a_t<b_t$.
\end{definition}
To compare $\omega$--vectors we allow trailing zeros to  equalize 
length. As usual, in  either ordering we write $L<M$ to mean $M>L$.

Both right and left orderings  are total and  
induce partial orderings  on  blocks by ordering their 
$\omega$--vectors.
We shall also talk about the left and right ordering of blocks. The 
\emph{reverse\/} lexicographic order is chosen in the right order case 
to provide consistency with the action of the Steenrod algebra, as we 
shall see later in \fullref{orderdown}. 

The following statement  is a simple numerical fact 
about unique descending $\omega$--vectors that will be required later.
\begin{proposition}\label{num}
If   ${\bf P}^d(n)$ admits a unique descending  $\omega$--vector  
$\omega$, then,  for any  block $B$ in  ${\bf P}^d(n)$  with   
$\omega(B)>_l\omega$, the first number $t$ for which   
$\omega_t(B)>\omega_t$ also satisfies the conditions $t>1$ and     
$\omega_{t-1}(B)<\omega_{t}(B)$.
\end{proposition}
For example,  ${\bf P}^8(3)$ has  the  unique $\omega$--vector 
$(2,1,1)$, with 
Ferrers block $F$,  and the left greater block $B$ as shown below.
\begin{example}
$$
F=\begin{smallmatrix}
1&1&1\\
1\\
0
\end{smallmatrix}\qua
B=\begin{smallmatrix}
1&1\\
1&1\\
0&1
\end{smallmatrix}
$$
We see that $\omega_1(B)=\omega_1(F)$ and $ \omega_2(B)>\omega_2(F)$. 
Also
$\omega_1(B)<\omega_2(B)$.
\end{example}
 
A familiar process in the combinatorics of Young tableaux 
is  \emph{straightening\/} which is a device for maneuvering a Young 
tableau into an equivalent   sum of  semistandard Young tableaux in 
the context of group rings.  We shall now explain an analogous 
process 
for blocks which we shall  later relate to the action of the Steenrod 
algebra on polynomials.   The idea is to  maneuver  a block in 
${\bf P}(n)$ into  a formal sum
 of semistandard blocks. We adopt the usual  notation $F_{i,j}$ for 
the $(i,j)$th entry of the block  $F$. 
   
\begin{definition}\label{ksplice}
Let $F$ be a block with $n$  rows  and let $k,t$ be  integers 
with $1\le k\le n$ and $t\ge 0$. Assume that, for a certain pair of
non-intersecting  sets  $S,T$, each containing   $k$   numbers 
between $1$
and $n$,  $F$ has entries $F_{i,t+2}=1$ and  $F_{i,t+1}=0$ for
$i\in S$ and  $F_{i,t+2}=0$ and  $F_{i,t+1}=1$ for $i\in T$. Let 
$G(S,T)$ be the
matrix formed from $F$ by leaving all entries
unchanged except in columns $t+1$ and $t+2$, where 
$G_{i,t+2}=0$ and $G_{i,t+1}=1$ for $i\in S$ and  $G_{i,t+2}=1$ and  
$G_{i,t+1}=0$ for $i\in T$. The process of replacing $F$ by the 
formal sum of the blocks
$G(S,T)$ for $S$ fixed and all possible $T$ is called  
\emph{$k$--splicing\/} of $F$ at \emph{column position\/}  $t+2$  and 
\emph{row positions\/} $S$.
\end{definition}
To put it briefly,  splicing replaces the block $F$ 
with the formal sum of  all the blocks $G(S,T)$  formed from 
$F$ by pulling  a selection of  $k$ digits  $1$ in column $t+2$ back  
one place into zero  positions, and pushing  a non-overlapping  
collection of $k$ digits in column $t+1$ forward one place  into zero 
 positions. Only two adjacent columns of the block  are altered, so 
effectively  splicing is a process carried out on a $2$--column block  
implanted as adjacent columns in a larger matrix.   Of course, even 
when the first part of the procedure is possible,  it may not always 
be possible to carry out the second part, in which case we define the 
result to be $0$. It should be noted that the $\alpha$   and  
$\omega$--vectors  of
each $G(S,T)$ are  the same as those of $F$. 

In the following example there is only one way of carrying out  
$2$--splicing of the matrix $B$ and the result is the sum of blocks 
$C$.
\begin{example}
$$B=\begin{smallmatrix}
0&1\\
1&0\\
0&1\\ 
1&0\\
1&0
\end{smallmatrix}
\qua
C=\begin{smallmatrix}
1&0\\
0&1\\
1&0\\ 
0&1\\
1&0
\end{smallmatrix}\qua +\qua
\begin{smallmatrix}
1&0\\
0&1\\
1&0\\ 
1&0\\
0&1
\end{smallmatrix}\qua +\qua
\begin{smallmatrix}
1&0\\
1&0\\
1&0\\ 
0&1\\
0&1
\end{smallmatrix}
$$
Here $S=\{1,3\}$ and there are three possible choices of $T$ 
corresponding to 
picking two rows from the list $\{2,4,5\}$.
\end{example}

We are now ready to establish the combinatorial part of our main 
theorem.
\begin{theorem}\label{mainbrace}
By  iterated  splicing, any  block with  descending $\omega$--vector 
can be 
replaced by a  formal  sum of  semistandard blocks.
\end{theorem}
\begin{proof}
Let $F$  be a block with $n$ rows   and descending $\omega$--vector. 
We argue by induction on  rows,  working from the bottom row upwards. 
We recall that for any block $F$ the subblock of the first $i$ rows 
of $F$
is denoted by $F[i]$.  As the inductive step,  assume that,  for 
some number $r$,  $F$ has been replaced by a formal sum of blocks $G$ 
such that the $\omega$--vectors  $\omega(G[i])$ of the subblocks are  
descending for  all
$i$ satisfying $1\le r+1\le i\le n$. 
The start of the induction is $r=n-1$,  since we are given that  
$\omega(F)$ is descending. If $r>0$ and  $\omega(G[r])$ fails to be 
descending,  then we can find a  column position $t+1$ for $t\ge 0$ 
such that 
$\omega_{t+1}(G[r])<\omega_{t+2}(G[r])$. Let $S$ denote the set of  
row positions 
$i$ in $G[r]$ for which $G[r]_{i,t+1}=0$ and  $G[r]_{i,t+2}=1$,  and 
suppose $S$ has $k$ elements. The effect of performing a $k$--splice 
of $G$ at column position $t+2$ and row positions $S$ is to produce a 
formal sum of matrices $H=H(S,T)$ with the properties
\begin{enumerate}
\item $H[i]$ has descending $\omega$--vector for $ r+1\le i\le n$,
\item $\omega(H[r]) >_l\omega(G[r])$.
\end{enumerate}
Assuming these two facts for the moment, we see by (i) that 
the process of splicing can be continued at column positions where 
the subblocks at level $r$ 
fail to be descending without disturbing the condition of descending 
$\omega$--vectors for levels below row $r$.   By (ii) this 
process must come to a stop since the $\omega$--vectors are bounded 
above in left  order. The process stops when all blocks $H$ are such 
that  $\omega(H[r])$  is descending   and this completes the 
inductive step. 

It remains to justify (i) and (ii). To obtain a 
typical block $H$ from $G$ we move $k$ digits $1$ of
$G[r]$ from column $t+2$ back to column $t+1$ and, say,  $a$  digits 
$1$ from column $t+1$ forward to column $t+2$. Since 
$\omega_{t+1}(G[r])<\omega_{t+2}(G[r])$ we must have $a<k$.  The 
other $k-a$ digits are moved from column $t+1$ to column $t+2$ below 
level $r$ in $G$.  Then 
$\omega_{t+1}(H[r])= \omega_{t+1}(G[r]) + k-a$.  It follows that  
$\omega_{t+1}(H[r])- \omega_{t+1}(G[r])>0$  which  proves (ii). 
Furthermore, let $i$ be a number between $r+1$ and $n$ and 
suppose  $b$ digits $1$, in  rows $r+1$ to $i$,  move from column 
$t+1$ to column $t+2$. Then   
$$ \omega_{t+1}(H[i])= \omega_{t+1}(G[i]) + k-a-b,\quad  
\omega_{t+2}(H[i])= \omega_{t+2}(G[i]) - k+a+b.$$
Hence 
$$ \omega_{t+1}(H[i])-  \omega_{t+2}(H[i])= \omega_{t+1}(G[i])- 
\omega_{t+2}(G[i]) +2( k-a-b).$$
By assumption we have  
$$a+b\le k \text{ and }\omega_{t+1}(G[i])\ge  
\omega_{t+2}(G[i]).$$ Hence $ \omega_{t+1}(H[i])\ge  
\omega_{t+2}(H[i])$ and since no other columns 
besides $t+1$  and $t+2$ have been disturbed, (i) follows.          
\end{proof}

\section{The hit problem for the Steenrod algebra}
\label{sec:sec3}
In this section we explain the action of the Steenrod algebra on 
polynomials and show how the combinatorial process of $k$--splicing   
can be realized by this action up to certain error terms. In 
favourable situations the error terms are  hit,  and this leads to 
our main theorem for  generators in row-regular degrees.
Background material on the hit problem can be found in 
Wood \cite{w1,w2,w3,w4,w5}, Janfada--Wood \cite{jw1,jw2},
Alghamdi--Crabb--Hubbuck \cite{ach} and Crabb--Hubbuck \cite{ch}.

The Steenrod algebra  ${\cal A}$ is a graded algebra generated by the 
Steenrod squares $\Sq^k$ in grading  $k$, over the field  ${\mathbb  
F}_2$,  subject to
the Adem relations (see Steenrod--Epstein \cite{s3}) and $\Sq^0=1$. 
\begin{proposition}\label{action}
The  Steenrod squares  $\Sq^k,k\ge 0 $, act    on  polynomials by 
linear  transformations  $\Sq^k\co  {\bf  P}^d\rightarrow 
 {\bf P}^{d+k}$, determined by the conditions,  
$$\Sq^1(x_i)=x_i^2,\qua \Sq^k(x_i)=0 \text{ for  }k>1,$$
 and  the Cartan formulae  for  polynomials $f,g$
$$\Sq^k(fg)=\sum_{i=0}^{k}\Sq^i(f)\Sq^{k-i}(g).$$
\end{proposition}
The  action of a general element of   ${\cal A}$ is by addition of 
compositions  of the Steenrod squares. The polynomial algebra  ${\bf 
P}(n)$ is a graded left  ${\cal A}$--module, where the grading is 
given by degree of polynomials. In principle  a Steenrod square can 
be evaluated on a monomial by iterated use of the Cartan formulae. A 
more compact way of stating
 these  formulae is in terms of the \emph{total squaring operation\/}, 
which is
 the formal sum  $SQ=1 +\Sq^1 + \Sq^2 + \cdots$. Then $SQ$ is 
multiplicative, that is, $SQ(fg)=SQ(f)SQ(g)$ for polynomials $f,g$ and the 
Cartan  formulae arise by comparing terms of degree $k$.

\begin{definition}
Two homogeneous polynomials $f,g$  of the same degree are 
\emph{equivalent  modulo hits\/} if they satisfy the relation
$$f=g +\sum_{i>0}\Sq^i (h_i), $$
over ${\mathbb F}_2$, 
which we refer to  as a \emph{hit equation\/}. In particular, if $g=0$  
then $f$
is hit. 
We write $f\cong g$  if $f-g$  is hit.
\end{definition}
The hit problem is to find a minimal generating set for the  
${\cal A}$--module ${\bf P}(n)$. Equivalently we want a vector space basis 
for the  quotient ${\bf Q}^d(n)$ of  $ {\bf  P}^d(n)$ by the hits in 
each degree $d$, frequently referred to as the \emph{cohits\/}.  Such a
 basis may  be represented by a list of monomials of degree $d$,  as 
in \fullref{ex02},  where  ${\bf Q}^4(3)$ has dimension $8$. 

The action of the Steenrod squares as described in \fullref{action},  when 
applied to polynomials  in an arbitrary number of variables, 
faithfully represents the Steenrod algebra in the sense that all 
relations in ${\cal A}$ can be detected by the action. Some 
elementary consequences for a  homogeneous polynomial
$f$  are easy to prove by induction on degree. 
\begin{proposition}\label{fractal}
If $k>\deg(f)$  then $\Sq^k(f)=0$,  and if  $k=\deg(f)$  then 
$\Sq^k(f)=f^2$. 
If $r$ is not divisible by $2^k$ then   $\Sq^r(f^{2^k})=0$
while   $\Sq^{s2^k}(f^{2^k})=(\Sq^s(f))^{2^k}$.
\end{proposition}
The second  statement expresses the \emph{fractal\/} nature of the 
Steenrod 
action. We shall  frequently invoke it  when considering the action 
of a Steenrod 
square on a monomial $b$ in terms of its columnwise action on the 
associated  block $B$.

There are some  important facts  about the Steenrod algebra which 
are  not immediately obvious from its  action  on polynomials. The 
Steenrod algebra ${\cal A}$ is multiplicatively generated by the 
Steenrod squares $\Sq^{2^k}$ for $k\ge 0$. It admits a coproduct which 
makes ${\cal A}$ into  a Hopf algebra (see Steenrod--Epstein \cite{s3}) with  a 
\emph{conjugation\/} operator  $\chi$. This is   a
 grade-preserving anti-automorphism of order 2. As in 
\fullref{action} there are rules for working out conjugates of  
Steenrod squares on
polynomials (see Walker--Wood \cite{ww}).
\begin{proposition}\label{chiaction}
The  action of the   conjugate  Steenrod squares  on  polynomials
\linebreak
$ \chi(\Sq^k)\co  {\bf  P}^d\rightarrow 
 {\bf P}^{d+k}$ are  determined by 
$$\chi(\Sq^k)(x_i)=x_i^{2^k}  \text{ if }k=2^a-1, a\ge 0,
\text{ and zero otherwise},
$$
and  the Cartan formula  $\chi(SQ)(fg)=\chi(SQ)(f)\chi(SQ)(g)$ for 
the  total conjugate square $\chi(SQ)=1 + \chi(\Sq^1) +  \chi(\Sq^2) 
+\cdots$.
\end{proposition}

There is one fact about the action of $\Sq^k$ and its conjugate 
$\chi(\Sq^k)$ on
a product of distinct variables that we shall need at a later stage 
in relation to the splicing process. Let $\{y_1,\ldots,y_m\}$ be a 
subset of the variables   $\{x_1,\ldots,x_n\}$.
\begin{lemma}\label{stsplice}
For $k\le m$,
$$\Sq^k(y_1\cdots y_m)=y_1\cdots y_m\sum_{\{i_1,\cdots,i_k\}}  
y_{i_1}\cdots y_{i_k},$$
where the summation is taken over $k$--element subsets of 
$\{y_1,\ldots,y_m\}$.
If $k>m$ then the result is zero. 
$$\chi(\Sq^k)(y_1\cdots y_m)= \Sq^k(y_1\cdots y_m) +f,$$
where every monomial in the polynomial $f$ has an exponent $\ge 4$.
\end{lemma}
\begin{proof}
By \fullref{chiaction} we have 
$$\chi(SQ)(y_1\cdots y_m)=\prod_{i=1}^m \chi(SQ)(y_i)=\prod_{i=1}^m 
(y_i+ y_i^2 + y_i^4 + \cdots).$$
Hence
$$\chi(SQ)(y_1\cdots y_m)=\prod_{i=1}^m 
(y_i+ y_i^2) + f  =SQ(y_1\cdots y_m) +f, $$
where all monomials in $f$ have  an exponent  $\ge 4$. The result 
then follows by 
comparing  terms of degree $m+k$.
\end{proof}

The following result has been significant in proving many
results on the hit problem and is known as the {\em $\chi$--trick} 
(Crossley \cite{c} and Wood \cite{w1,w2,w3,w4,w5}).
\begin{proposition}\label{chi}
For homogeneous polynomials $u,v$ 
$$u\Sq^k(v) -v\chi(\Sq^k)(u)  = \sum_{i>0} \Sq^i(v \chi(\Sq^{k-i})(u)).$$
\end{proposition}
Thus  $u\Sq^k(v)\cong v\chi(\Sq^{k})(u)$ and the statement extends by 
composition and addition of Steenrod operations to show that 
$u\Theta(v)\cong v\chi(\Theta)(u)$  for any element $\Theta$ in $\cal 
A$. The $\chi$--trick is the analogue of integration by parts in 
calculus, when the Steenrod squares are interpreted as differential 
operators (see Wood \cite{w6}).  An immediate application of the $\chi$--trick is 
the following well known  observation,
used  to prove  the Peterson conjecture \cite{w1,w2,w3,w4,w5}.

\begin{proposition}\label{excess}
Let $u$ and $v$ be homogeneous polynomials such that $\deg(u)<\mu(\deg(v))$. 
Then 
$uv^2$ is hit. 
\end{proposition} 
The proof follows by writing $v^2=\Sq^d(v)$, where $d$ is the degree 
of $v$, applying the $\chi$--trick,  and then the fact that the excess 
 of  $\chi(\Sq^d)$ is $\mu(d)$. The condition  $\deg(u)<\mu(d)$ and 
the definition of excess (see Steenrod--Epstein \cite{s3}) implies  $\chi(\Sq^d)(u)=0$.
 
We shall find it convenient  to 
switch back and forth between blocks and  monomials where 
appropriate. To avoid repetition we adopt the temporary convention of 
using  upper case letters for  blocks  and their lower case versions  
for corresponding monomials. 
A  vertical partition of a block
$B=FG$ corresponds to the monomial $b=fg^{2^t}$ if $F$ has $t$ 
columns. If $t=0$ then $F$ is empty (corresponding to the monomial 
$1$).  If $H=\sum H_k$ is a formal sum of blocks then $FH$ is the 
formal sum  $\sum FH_k$. From the Cartan formula  of  
\fullref{action} and the fractal nature of the action of Steenrod 
squares,  as explained in   \fullref{fractal},  we can  study the 
action of $\Sq^m$ columnwise on blocks. For example, corresponding to 
$B=FG$ we have 
$$\Sq^{m}(b) = \sum \Sq^p(f)(\Sq^q(g))^{2^t},$$
where the summation is over all $p,q\ge 0$ with $p+2^tq =m$. 
Splitting a block into its columns  as $B=B_1B_2\ldots B_t$  leads to 
the formula
$$\Sq^{m}(b) = \sum \Sq^{p_1}(b_1)(\Sq^{p_2}(b_2))^2\ldots 
(\Sq^{p_t}(b_t))^{2^t},$$
where the summation is taken over all solutions in non-negative 
integers $p_i$
of the equation $p_1+2p_2 + \cdots +2^tp_t=m$.
In the light of \fullref{stsplice},  describing the action of a 
Steenrod square  on a product of distinct variables,  it is easy to 
see that the typical action  of a  Steenrod square on a block  moves  
 digits  $1$ from  one  column to the next  column on the right in 
the same row, with the knock-on effect of binary addition if digits 
superimpose. In particular we deduce the following fact about the 
order relations introduced in \fullref{sec:sec1}.
\begin{lemma}\label{orderdown}
Let $\Sq^k(F)=F_1+\cdots +F_s$, for $k\ge 1$, be a formal sum of 
distinct blocks.  Then  $F_i< F$  in both  the left  and right  
orderings for $1\le i \le s$.
\end{lemma}

We shall now interpret \fullref{stsplice} in block language and use 
the $\chi$--trick to show how $k$--splicing in the second column of a 
$2$--column block can be  realized by the action of the  Steenrod 
algebra modulo certain error terms.  

Let $C$ be a $2$--column block. Let   $ R$ be the set  of row 
positions   where there is a digit $1$ in $C$.  Partition $ R$  into 
three
subsets  ${ U,V,W}$ as follows.  For $i\in  { U}$ we require 
$C_{i,1}=C_{i,2}=1$, and for   $i\in  { V}$ we require $C_{i,1}=0$  
and $C_{i,2}=1$,  and for   $i\in  { W}$ we require $C_{i,1}=1$ and 
$C_{i,2}=0$. Now select a subset ${ S}$ of $k$ elements of ${ V}$ and 
let $B$ be the $2$--column  block with zero entries except for 
$B_{i,1} =1$ for $i\in {S}$. Let $A$ be the block formed from $C$ by 
deleting the  digits in positions $C_{i,2}$ for $i\in { S}$. The 
following diagrams illustrate an example where
$$
{U}=\{1\}, \qua {V}=\{2,3,4\},
\qua  {W}=\{5,6,7\},\qua  { S}=\{2,3\},\qua k=2.$$
\begin{example}\label{ex3}
$$
C=
\begin{smallmatrix}
1&1\\
0&1\\
0&1\\
0&1\\ 
1&0\\
1&0\\
1&0
\end{smallmatrix}\qua
A=
\begin{smallmatrix}
1&1\\
0&0\\
0&0\\
0&1\\ 
1&0\\
1&0\\
1&0
\end{smallmatrix}\qua
B=
\begin{smallmatrix}
0&0\\
1&0\\
1&0\\
0&0\\ 
0&0\\
0&0\\
0&0
\end{smallmatrix}
$$
\end{example}

Then $c=a\Sq^k(b) \cong b\chi(\Sq^k)(a)$ by the $\chi$--trick.  Now 
monomials with exponents $\ge 4$ correspond to blocks which are right 
lower than any  $2$--column block. Hence by   \fullref{stsplice} the 
effect of  $\chi(\Sq^k)$ on $A$ is the same as $\Sq^k$ on $A$  modulo 
right lower blocks.
Furthermore,  any  effect arising from  $\Sq^k$ 
via the Cartan formula on the second column of $A$ also produces 
right lower  blocks, as does the action of $\Sq^k$ on the first column 
on any row in the set
 ${ U}$ by the knock-on effect of binary addition. On the other hand, 
the effect of   $\Sq^k$ concentrated on the rows of the set ${ W}$ is 
to  produce exactly the result  of the $k$--splicing process on $A$.  
Consequently 
$b\chi(\Sq^k)(a)$ produces the effect of $k$--splicing $C$ modulo the 
error terms as described. This is summarized in the following 
statement. 
\begin{proposition}\label{bstsplice}
Let $C$ be a $2$--column block and let $C'$ be the sum of $2$--column 
blocks arising from a $k$--splicing process of $C$ at  column $2$. 
Then $C\cong C'$ modulo blocks which are right lower  than any 
$2$--column block.
\end{proposition}

Now we need to investigate what happens when a block is implanted as 
adjacent columns in a larger  block. 
\begin{lemma}\label{midorder}
Let $B=FCG$ be a vertical splitting of a block and suppose $C\cong C' 
+ R$, where $C'$ is a sum of blocks of the same size as $C$,  and $R$ 
is a sum of blocks each of which is  right lower than $C$.  Let 
$B'=FC'G$. 
Then $B\cong B' + F'H +FK$,
where $F'$ is a sum of blocks of the same size as $F$, each of which 
is left lower than $F$,  and $K$ is a sum of blocks each of which is 
right lower than
$CG$. In particular $B$ is equivalent to $B'$ modulo blocks which 
are either left or right lower than $B$. 
\end{lemma}
\begin{proof}
Substituting $R$ for $C$ in $B$ immediately produces blocks which may 
overlap with $G$ but certainly have the form $FK$,
as stated in the proposition. We  may therefore  assume that  $R=0$.
In terms of corresponding monomials we have $b=fc^{2^t}g^{2^s}$, 
where $t$ is the number of columns in $F$ and $s-t$ the number of 
columns in $C$.  Then  $C\cong C'$ and  there is a hit equation $c=c' 
+\sum_{k>0} \Sq^kh_k$. By the fractal property in \fullref{fractal}, 
we have the hit equation
$c^{2^t}=(c')^{2^t} + \sum_{k>0} \Sq^{2^tk}(h_k^{2^t})$. Applying the 
$\chi$--trick in \fullref{chi} to $u=fg^{2^s}$ and $v=h_k^{2^t}$ for 
each $k$ in turn and then adding, we see that $b-b'\cong \Theta(u)v$ 
for some positively graded element $\Theta$ in the Steenrod algebra. 
Then by the Cartan formula,  $\Theta$ must have a positive action 
either on $f$ or $ g^{2^s}$ which means that, in the language of 
blocks,  by \fullref{orderdown},  either  $F$ or $G$ is moved to a 
sum of lower blocks in either order. The result follows.
\end{proof}
An immediate corollary of \fullref{midorder} and \fullref{bstsplice} 
is the following result.

\begin{proposition}\label{link}
Let $B=FCG$ be a partitioned block, where $C$ has two columns in 
positions
$t+1,t+2$. Let $C'$ be the sum of the $2$--column  blocks arising out 
of a  $k$--splicing process 
at  column $t+2$ and let $B'=FC'G$.  Then  $B\cong B'$ modulo blocks 
of the form $F'H$, where $F'$ is left lower than $F$, and blocks  
$FK$, where $K$ is right lower than $CG$. In particular $B\cong B'$ 
modulo blocks  which are either left or right lower than $B$. 
\end{proposition}
Of course, if for some choice of $k$, the first stage of $k$--splicing 
is not possible, then the above statement is void. On the other hand, 
if $k$ can be chosen in such a way that  $k$--splicing produces the 
zero result,
then  \fullref{link} says that $B$ is reducible modulo hits to blocks 
which are either left or right lower than $B$ in the specified way. 
This leads to the following result.

\begin{proposition}\label{lessleft}  
Let  $B$ in ${\bf P}^k(n)$ be a block whose $\omega$--vector is not 
descending, 
so that  $\omega_{t+1}(B)<\omega_{t+2}(B)$ for some value of $t$.  We 
can write $B=FCG$, where $F$ has $t$ columns and $C$ is a $2$--column 
block with $\omega_1(C)<\omega_{2}(C)$. Then  $B$ is hit modulo  
blocks of the form $F'H$, where $F'$ is left lower than $F$, and 
blocks  $FK$, where $K$ is right lower than $CG$. In particular $B$ 
is hit modulo blocks  which are either left or right lower than $B$.  
\end{proposition}
\begin{proof}
The condition  $\omega_1(C)<\omega_{2}(C)$ ensures  $k$--splicing of 
$C$ in the second column is possible and the largest such $k$ 
produces the zero result. 
\end{proof}

We shall now exploit the above  results  in a situation where we can
control the error terms. Parts of the next proposition, originating in
Singer's work \cite{s4}, are known in more generality (see Carlisle--Wood
\cite{cw} and Mothebe \cite{m2}) but for the sake of completeness we
include proofs of these  particular cases.

\begin{proposition}\label{lowlspike}
Assume that  ${\bf P}^d(n)$ admits a unique  descending 
$\omega$--vector $\omega$. Let $B$ be a block in ${\bf P}^d(n)$.  
\begin{enumerate}
\item If $\omega(B)<_l\omega $ then $B$ is hit. 
\item If $\omega(B)<_r\omega $ then $B$ is hit. 
\item There is a generating set of blocks $B$  for  
${\bf Q}^d(n)$ with $\omega(B)=\omega$.
\end{enumerate}
\end{proposition}

\begin{proof}
We start with the proof of (i).  Let $B$ be a block in  
${\bf P}^d(n)$ with  $\omega(B)<_l\omega$. There is a first column 
position $t<n$ from the left  where   $\omega_t(B)<\omega_t$. 
Consider the vertical splitting $B=FCG$ where $F$ has $t-1$ columns 
(empty if $t=1$), $C$  has one column in position $t$,  and $G$ 
(non-empty) has the rest of the columns of $B$.
Then $\mu(\deg(G))>\deg(C)$, otherwise  we can create a block 
$B'=FCG'$, where $G'$ is a spike with fewer rows than   $C$. Then  
$\omega(B')$ is descending and $\omega(B')<_l\omega$,  contrary to 
the assumption that $\omega$ is the 
unique  descending  $\omega$--vector. It follows from \fullref{excess}
that $CG$ is hit. Then by the  arguments used in previous 
propositions we see 
that $B\cong F'H$, where $F'$ is a sum of blocks of the same size as 
$F$ and lower than $F$  in the left order. For a  typical such block 
$B''$ we have    $\omega_s(B'')<\omega_s$ for some $s<t$. Iteration 
of the process must come to a 
stop at or before $t=1$ when the result is zero. Hence $B$ is hit.  

To prove (ii) we may as well start with a block $B$ for which 
$\omega(B)<_r\omega $ and $\omega(B) >_l \omega $. Let $t$ be the 
first
number such that  $\omega_{t+2}(B) > \omega_{t+2} $. Then by 
\fullref{num}
we can write $B=FCG$ where $F$ (possibly empty) has $t$ columns and 
$C$ has 
two columns with  $\omega_1(C)<\omega_{2}(C)$ as in 
\fullref{lessleft}.
According to this proposition $B$ is   equivalent to a sum of  blocks 
of the form $F'H$, where $F'$ is left lower than $F$, and blocks  
$FK$, where $K$ is right lower than $CG$. But then  $\omega(F'H)$ is 
left lower than $\omega $ and therefore $F'H$ is hit by part (i). 
Hence  $B$ is equivalent to a sum of blocks right lower than 
$B$. In particular their $\omega$--vectors are right lower than 
$\omega$ and the process can therefore be iterated.  The  procedure 
must come to a stop since we cannot have an infinite chain of   right 
lower blocks. The process ends when the result is zero,   and this  
proves that $B$ is hit. 

The proof of (iii) follows the same line of argument  as the 
proof of (ii),  except that the process stops when the 
$\omega$--vectors of the blocks reach 
$\omega$. 
\end{proof}
We now state and prove the main result.
\begin{theorem}\label{main}
Suppose that ${\bf P}^d(n)$ admits a unique descending 
$\omega$--vector. 
Then the cohits  ${\bf Q}^d(n)$ are spanned by the semistandard 
blocks.   
\end{theorem}
\begin{proof}
By part (iii) of \fullref{lowlspike} we can start with a 
spanning set for   ${\bf Q}^d(n)$ consisting of blocks $B$ having the 
unique descending 
$\omega$--vector. By \fullref{link} and \fullref{lowlspike} we 
can replace  $B$ by the result of any $k$--splicing modulo hits. The 
proof is then complete by \fullref{mainbrace}. 
\end{proof}

The limitation  of the above approach in the non-regular case,  where 
there 
is more than one descending $\omega$--vector,  is  illustrated by the 
 example   ${\bf P}^7(4)$. Here  there are two descending   
$\omega$--vectors $(1,1,1)<(3,2)$, the least and greatest in either of 
the order relations.  The other possible  $\omega$--vectors are 
$(3,0,1)$ and $(1,3)$ which lie between these extremes. 
\begin{example}  
Consider the following  block $C$ with  $\omega(C)= (1,3)$.  
$$
C=
\begin{smallmatrix}
1\\
0&1\\ 
0&1\\
0&1
\end{smallmatrix}\qua
E=
\begin{smallmatrix}
0&0&1\\
1\\ 
1\\
1
\end{smallmatrix}\qua
 F=
\begin{smallmatrix}
0&1\\
0&1\\ 
0&1\\
1
\end{smallmatrix}\qua
G=
\begin{smallmatrix}
0&1\\
0&1\\ 
1\\
0&1
\end{smallmatrix}\qua
H=
\begin{smallmatrix}
0&1\\
1\\ 
0&1\\
0&1
\end{smallmatrix}
$$ 
Now $3$--splicing of $C$ in the second column  has zero effect but the 
Steenrod realization has error term $E$ with  $\omega(E)= (3,0,1)$. 
Hence  $C\cong E$.
Similarly,  $1$--splicing $E$ in the third column produces the 
equivalence
$E\cong F+G+H$. So iterated splicing has produced the relation
$$C+F+G+H\cong 0,$$
involving blocks with  $\omega$--vector $(1,3)$.
However, it can  be shown  that $C$ is not equivalent to a 
combination of blocks with    $\omega$--vectors $(1,1,1)$ or $(3,2)$.
\end{example}
This example contrasts with the case $n=3$, where a basis for the 
cohits
 can  be taken with descending   $\omega$--vectors. The complete 
solution of 
the hit problem in the case $n=4$ has been given by Kameko \cite{k3} 
in a 
format which analyzes the hit problem one $\omega$--vector at a time. 
The vector space ${\bf Q}^{\omega}(n)$ is formed  by taking the 
quotient of the subspace of   ${\bf P}^{d}(n)$ generated by monomials 
with $\omega$--vector $\le \omega$ by the hits and the subspace 
generated by monomials with  $\omega$--vector  $<\omega$ in left 
order.  Much of the above work can be applied   to  ${\bf 
Q}^{\omega}(n)$ 
when  $\omega$ is the  least descending  $\omega$--vector in degree 
$d$ (which is the same in either order).  

\section{The Steinberg representation}
\label{sec:sec4}
The  degree $d=2^n-n-1$ is row-regular for $n$  and
 $d=\sum_{i=1}^n (2^{n-i}-1)$ is the unique exponential partition of 
$d$ into  $n$ parts, with unique descending $\omega$--vector 
$(n-1,n-2,\ldots,1,0)$ and 
Ferrers block  $F$. The corresponding partition is also 
$\lambda=(n-1,n-2,\ldots,1,0)$. The number of semistandard Young 
tableaux, and therefore semistandard blocks,  is given in Fulton \cite[page 
55]{f} by the \emph{hook formula\/} 
$$d_{\lambda}(m)=\prod_{(i, j)\in \lambda} \frac{m+j-i}{h(i,j)},$$
for the Ferrers diagram of     $\lambda$, filled with numbers from 
the set $\{1,\ldots,m\}$, 
where  $h(i,j)$ denotes the hook length of the node  in the Ferrers 
diagram at 
position $(i, j)$, that is,  the number of nodes to the right and below 
the given position in the Ferrers diagram including the position 
itself.

In our application,  $m=n$ and $h(i,j)=2(n-i-j)+1$ giving
$d_{\lambda}(n)=2^{\binom{n}{2}}$, the dimension of the Steinberg
representation of  $GL(n,{\mathbb F}_2)$ (see Mitchell--Priddy
\cite{mp}). \fullref{main} shows that the dimension of the vector space of
cohits   ${\bf Q}^d(n)$ is bounded by  $2^{\binom{n}{2}}$. The remarks
in \fullref{sec:sec1} about the first occurrence of an irreducible
representation then  finally establish \fullref{th01}.

For $m<n$,  the Weyl module for  $GL(n,{\mathbb F}_2)$ corresponding 
to the partition  $\lambda=(m-1,m-2,\ldots,1,0,\ldots,0)$ is 
irreducible, and has dimension   $d_{\lambda}(n)$. By 
Carlisle--Kuhn \cite[Theorem 1.1]{ck}, the first occurrence 
as  a composition factor is in degree  $d=2^{m+1}-1-m$. The work 
above can then be applied to  ${\bf Q}^{\omega}(n)$,  
when   $\omega$ is the least descending  $\omega$--vector in degree 
$d$, to show that 
$\dim({\bf Q}^{\omega}(n)) = d_{\lambda}(n)$.  

\bibliographystyle{gtart}
\bibliography{link}

\begin{thebibliography}{}
\providecommand\bibmarginpar{\leavevmode\marginpar}
\def\urlstyle#1{{\tt #1}}

\bibitem{ach}
\textbf{M\,A Alghamdi}, \textbf{M\,C Crabb}, \textbf{J\,R Hubbuck},
  \emph{Representations of the homology of {$BV$} and the {S}teenrod algebra
  I}, from: ``Adams Memorial Symposium on Algebraic Topology 2 (Manchester,
  1990)'', London Math. Soc. Lecture Note Ser. 176, Cambridge Univ. Press,
  Cambridge (1992)  217--234 \xox{MR}{1232208}

\bibitem{b}
\textbf{J\,M Boardman}, \emph{Modular representations on the homology of powers
  of real projective space}, from: ``Algebraic topology (Oaxtepec, 1991)'',
  Contemp. Math. 146, Amer. Math. Soc., Providence, RI (1993)  49--70
  \xox{MR}{1224907}

\bibitem{ck}
\textbf{D Carlisle}, \textbf{N\,J Kuhn},
  \href{http://dx.doi.org/10.1016/0021-8693(89)90073-2} {\emph{Subalgebras of
  the {S}teenrod algebra and the action of matrices on truncated polynomial
  algebras}}, J. Algebra 121 (1989) 370--387 \xox{MR}{992772}

\bibitem{cw}
\textbf{D\,P Carlisle}, \textbf{R\,M\,W Wood}, \emph{The boundedness conjecture
  for the action of the {S}teenrod algebra on polynomials}, from: ``Adams
  Memorial Symposium on Algebraic Topology, 2 (Manchester, 1990)'', London
  Math. Soc. Lecture Note Ser. 176, Cambridge Univ. Press, Cambridge (1992)
  203--216 \xox{MR}{1232207}

\bibitem{ch}
\textbf{M\,C Crabb}, \textbf{J\,R Hubbuck}, \emph{Representations of the
  homology of {$BV$} and the {S}teenrod algebra. {II}}, from: ``Algebraic
  topology: new trends in localization and periodicity (Sant Feliu de Gu\'\i
  xols, 1994)'', Progr. Math. 136, Birkh\"auser, Basel (1996)  143--154
  \xox{MR}{1397726}

\bibitem{c}
\textbf{M\,D Crossley}, \emph{{$H^{*}V$} is of bounded type over
  {$\mathcal{A}_p$}}, from: ``Group representations: cohomology, group actions
  and topology (Seattle, WA, 1996)'', Proc. Sympos. Pure Math. 63, Amer. Math.
  Soc., Providence, RI (1998)  183--190 \xox{MR}{1603151}

\bibitem{f}
\textbf{W Fulton}, \emph{Young tableaux}, London Mathematical Society Student
  Texts 35, Cambridge University Press, Cambridge (1997) \xox{MR}{1464693}

\bibitem{j}
\textbf{G James}, \textbf{A Kerber}, \emph{The representation theory of the
  symmetric group}, Encyclopedia of Mathematics and its Applications 16,
  Addison-Wesley Publishing Co., Reading, MA (1981) \xox{MR}{644144}

\bibitem{jw1}
\textbf{A\,S Janfada}, \textbf{R\,M\,W Wood},
  \href{http://dx.doi.org/10.1017/S0305004102006059} {\emph{The hit problem for
  symmetric polynomials over the {S}teenrod algebra}}, Math. Proc. Cambridge
  Philos. Soc. 133 (2002) 295--303 \xox{MR}{1912402}

\bibitem{jw2}
\textbf{A\,S Janfada}, \textbf{R\,M\,W Wood},
  \href{http://dx.doi.org/10.1017/S0305004103006662} {\emph{Generating
  {$H^{*}({\rm BO}(3),\mathbb F\sb 2)$} as a module over the {S}teenrod
  algebra}}, Math. Proc. Cambridge Philos. Soc. 134 (2003) 239--258
  \xox{MR}{1972137}

\bibitem{k3}
\textbf{M Kameko}, \emph{Generators of the cohomology of $BV_4$}, preprint

\bibitem{k1}
\textbf{M Kameko}, \emph{Products of projective spaces as Steenrod modules},
  PhD thesis, Johns Hopkins University (1990)

\bibitem{k2}
\textbf{M Kameko}, \emph{Generators of the cohomology of {$BV_3$}}, J. Math.
  Kyoto Univ. 38 (1998) 587--593 \xox{MR}{1661173}

\bibitem{m}
\textbf{I\,G Macdonald}, \emph{Symmetric functions and {H}all polynomials},
  second edition, Oxford Mathematical Monographs, Oxford University Press,
  Oxford (1995) \xox{MR}{1354144}

\bibitem{mt}
\textbf{P\,A Minh}, \textbf{T\,T Tri},
  \href{http://dx.doi.org/10.1090/S0002-9939-99-05424-6} {\emph{The first
  occurrence for the irreducible modules of general linear groups in the
  polynomial algebra}}, Proc. Amer. Math. Soc. 128 (2000) 401--405
  \xox{MR}{1676308}

\bibitem{mp}
\textbf{S\,A Mitchell}, \textbf{S\,B Priddy},
  \href{http://dx.doi.org/10.1016/0040-9383(83)90014-9} {\emph{Stable
  splittings derived from the {S}teinberg module}}, Topology 22 (1983) 285--298
  \xox{MR}{710102}

\bibitem{m2}
\textbf{M\,F Mothebe}, \href{http://dx.doi.org/10.1081/AGB-120003466}
  {\emph{Generators of the polynomial algebra {$F\sb 2[x\sb 1,\dots,x\sb n]$}
  as a module over the {S}teenrod algebra}}, Comm. Algebra 30 (2002) 2213--2228
  \xox{MR}{1904635}

\bibitem{p}
\textbf{F\,P Peterson}, \emph{Generators of
  $\mathbf{H}^{*}(\mathbb{R}P^{\infty}\wedge \mathbb{R}P^{\infty})$ as a module
  over the Steenrod algebra}, Abstracts Amer. Math. Soc. 833-55-89 (1987)

\bibitem{s1}
\textbf{B\,E Sagan}, \emph{The symmetric group}, second edition, Graduate Texts
  in Mathematics 203, Springer, New York (2001) \xox{MR}{1824028}

\bibitem{s4}
\textbf{W\,M Singer},
  \href{http://links.jstor.org/sici?sici=0002-9939(199102)111:2%3C577:OTAOSS%3%
E2.0.CO%3B2-K} {\emph{On the action of {S}teenrod squares on polynomial
  algebras}}, Proc. Amer. Math. Soc. 111 (1991) 577--583 \xox{MR}{1045150}

\bibitem{s2}
\textbf{R\,P Stanley}, \emph{Enumerative combinatorics. {V}ol. 2}, Cambridge
  Studies in Advanced Mathematics 62, Cambridge University Press, Cambridge
  (1999) \xox{MR}{1676282}

\bibitem{s3}
\textbf{N\,E Steenrod}, \textbf{D\,B\,A Epstein}, \emph{Cohomology operations},
  Annals of Mathematics Studies 50, Princeton University Press, Princeton, N.J.
  (1962) \xox{MR}{0145525}

\bibitem{ww}
\textbf{G Walker}, \textbf{R\,M\,W Wood},
  \href{http://dx.doi.org/10.1006/jabr.2001.8982} {\emph{Linking first
  occurrence polynomials over {$\mathbb F\sb 2$} by {S}teenrod operations}}, J.
  Algebra 246 (2001) 739--760 \xox{MR}{1872123}

\bibitem{w1}
\textbf{R\,M\,W Wood}, \emph{Steenrod squares of polynomials and the {P}eterson
  conjecture}, Math. Proc. Cambridge Philos. Soc. 105 (1989) 307--309
  \xox{MR}{974986}

\bibitem{w6}
\textbf{R\,M\,W Wood}, \href{http://dx.doi.org/10.1112/S0024611597000324}
  {\emph{Differential operators and the {S}teenrod algebra}}, Proc. London
  Math. Soc. $(3)$ 75 (1997) 194--220 \xox{MR}{1444319}

\bibitem{w2}
\textbf{R\,M\,W Wood}, \href{http://dx.doi.org/10.1112/S002460939800486X}
  {\emph{Problems in the {S}teenrod algebra}}, Bull. London Math. Soc. 30
  (1998) 449--517 \xox{MR}{1643834}

\bibitem{w3}
\textbf{R\,M\,W Wood}, \emph{Hit problems and the {S}teenrod algebra}, from:
  ``Proceedings of the summer school `Interactions between algebraic topology
  and invariant theory', a satellite conference of the third European congress
  of mathematics, Ioannina University, Greece'' (2000)  65--103

\bibitem{w4}
\textbf{R\,M\,W Wood}, \emph{Invariants of linear groups as modules over the
  {S}teenrod algebra}, from: ``Ingo2003, Invariant theory and its interactions
  with related fields, G\"ottingen'' (2003)

\bibitem{w5}
\textbf{R\,M\,W Wood}, \emph{The {P}eterson conjecture for algebras of
  invariants}, from: ``Invariant theory in all characteristics'', CRM Proc.
  Lecture Notes 35, Amer. Math. Soc., Providence, RI (2004)  275--280
  \xox{MR}{2066475}

\end{thebibliography}
\end{document}